\documentclass[titlepage,twoside,12pt]{article}
\usepackage{amsmath}
\usepackage{amssymb}
\usepackage{amsfonts}
\begin{document}

{\LARGE \bf Inevitable Infinite Branching in the \\ \\ Multiplication of Singularities} \\

{\bf Elem\'{e}r E ~Rosinger} \\ \\
{\small \it Department of Mathematics \\ and Applied Mathematics} \\
{\small \it University of Pretoria} \\
{\small \it Pretoria} \\
{\small \it 0002 South Africa} \\
{\small \it eerosinger@hotmail.com} \\ \\

\hfill {\it Dedicated to Marie-Louise Nykamp} \\ \\

{\bf Abstract} \\

Singularities appear in numerous important mathematical models used
in Physics. And in most of such cases singularities are involved in
essentially nonlinear contexts. For more than four decades, general
enough nonlinear theories of singularities have been developed. A
critically important related feature is that, above certain levels
in singularities, the operation of multiplication, and in general,
nonlinear operations on such singularities do inevitably branch in
infinitely many ways, without the possibility for the existence of
some unique natural or canonical way such nonlinear operations may
be performed. Consequently, the choice in such
branchings has to come from extraneous considerations. \\ \\

\hspace*{3cm} ''We do not possess any method at all to derive \\
\hspace*{3cm} systematically solutions that are free of \\
\hspace*{3cm} singularities...'' \\
\hspace*{3cm} Albert Einstein : The Meaning of Relativity. \\
\hspace*{3cm} Princeton Univ. Press, 1956, p. 165 \\ \\

{\bf 0. Preliminaries} \\

Singularities have been present in Mathematics ever since the simple and natural looking issue of dividing by zero. And with the mathematisation of modern Physics, they are causing major difficulties in quite a number of disciplines of that field of science. Needless to say, various branches of engineering, as well as other fields of science and technology encounter similar difficulties due to singularities in the mathematical models employed. \\

Here we shall consider singularities of functions $f : X \longrightarrow Y$, where $X$ and $Y$ are some Euclidean domains or finite dimensional manifolds. Also, we shall consider singularities of various generalized functions which extend the above concept of function. In this way, the following {\it two fundamental features} of singularities will be of concern, namely

\begin{itemize}

\item the extent of the set of singularities as a subset in the domain of definition of functions or generalized functions,

\item the behaviour of functions or generalized functions in the neighbourhood of singularities.

\end{itemize}

The main conclusion obtained will be as follows :

\begin{itemize}

\item In case no limitations are imposed on the above two features of singularities, as soon as multiplication, and in general, nonlinear operations are effectuated with generalized functions, there is an {\it inevitable infinite branching} in the way such operations can be defined. In other words, there is {\it no} canonical, natural or unique way multiplication, and in general, nonlinear operations can be defined for such singular generalized functions. Consequently, the specific choice of the result of multiplication, and in general, of nonlinear operations on generalized functions with singularities has to be made based on {\it extraneous} considerations.

\end{itemize}

A first systematic and far reaching mathematical approach to singularities was given in the 19th century by the theory of Functions of One Complex Variable. That was the time which led, among others and not necessarily in a related manner, to the celebrated and yet unsolved Riemann Hypothesis, which shows the depth of the respective theory. \\
As for singularities in the context of functions of one complex variable, one should recall the Great Picard Theorem, according to which an analytic function in the neighbourhood of an isolated singularity point that is an essential singularity will take on all possible complex values infinitely often, with at most a single exception. \\

Consequently, in the neighbourhood even of one single and isolated singularity, one can expect a rather arbitrary behaviour, when one deals with more general functions than analytic ones. Not to mention the situation when the singularity points form a considerably larger subset in the domain of definition of a function. Therefore, the consideration of the above two aspects related to singularities is indeed appropriate, and in fact, necessary in the case of a deeper going and more wide ranging approach. \\

Beyond the confines of analytic functions, it was the linear theory of Schwartz distributions which, starting with the mid 20th century, gave in certain respects a considerably more powerful and general treatment of singularities. \\

The major limitation of the Schwartz approach is in its essential confinement to linear situations, and its consequent inability to deal with singularities in a nonlinear context, and do so in a convenient and systemtic manner, without the recourse to what often are merely adhoc approaches. This contrast sharply with the remarkable natural ease and clarity the earlier complex function theory was applicable within a nonlinear context, restricting itself, as it did naturally, to analytic situations. \\
Furthermore, as a rather unfortunate event, the more notorious than celebrated 1954 paper of Schwartz [1] claimed to prove that a nonlinear theory of distributions would altogether be undesirable, if not in fact, impossible ... \\
Amusingly, this claim has attained a wider acceptance, [2], and consequently has for long distorted the perception of the situation concerning singularities within a general enough nonlinear context ... \\

As it happened, however, two general nonlinear approaches to
singularities have nevertheless emerged, the first being the one
listed under 46F30 in the Mathematical Subject Classification of the
American Mathematical Society, and with the second the Abstract
Differential Geometry, or ADG, established by A Mallios in the
1990s, see [3-11] for some of the more relevant references in this regard. \\

History, nevertheless, still seems to hang somewhat uneasily upon the issue of singularities ... \\
Indeed, due to the remarkable generality, clarity and power of the Schwartz linear distributional approach, the issue of algebraic operations on singularities has naturally and inevitably been {\it restricted} to addition alone, without a similarly general consideration of the operation of {\it multiplication}, an operation which ended up being in fact implicitly excommunicated form any possible suitable and general enough theory of singularities, in view of the mentioned 1954 Schwartz paper and its long ongoing misinterpretations. \\

And the {\it surprising} and hardly yet noted fact here is as follows. Although multiplication is closely related to addition, being in certain ways but a repeated addition, when it nevertheless comes to singularities, an {\it essential} difference appears between these two basic algebraic operations. Namely

\begin{itemize}

\item addition can be extended from usual functions, be they regular or with singularities, to all sort of generalized functions, and such an extension appears to be naturally done in a unique, canonical manner,

\end{itemize}

while on the other hand

\begin{itemize}

\item due to most simple algebraic, more precisely, ring theoretic reasons, the extension of multiplication from usual functions, be they regular or with singularities, to generalized functions does {\it no} longer have such a naturally unique canonical way.

\end{itemize}

Here is, therefore, the root of what is in fact no less than an {\it infinite branching} of the ways multiplications can naturally be defined for generalized functions. \\
And no wonder, this root has so far been mostly missed due to the mentioned implicit omission to consider multiplication of singularities within a general enough, and not merely adhoc, context ... \\

Within the general nonlinear treatment of singularities in 46F30, so far three classes of differential algebras of
generalized functions have been used in a variety of problems, mainly for the solution of large classes of nonlinear
systems of partial differential equations. \\
It is instructive to recall the way these three classes of algebras relate to the above two fundamental features of
singularities. \\
In this regard, all these algebras are able to deal with the singularities the Schwartz linear theory of distributions
can, since each of these algebras contains all the Schwartz distributions. \\

The issue, therefore, is to what extent these algebras are able to
deal with additional singularities. Let us then consider these
algebras defined on any given Euclidean open set $X \subseteq
\mathbb{R}^n$. Consequently, the generalized functions will extend
various classes of usual functions $f : X
\longrightarrow E$, where $E$ is a finite dimensional real or complex vector space. \\

Here it should be pointed out again that the class of admissible
singularities is large in no less than {\it two} significantly
useful ways :

\begin{itemize}

\item First, the singularities of the functions $f : X \longrightarrow
E$ considered can be given by arbitrary subsets $\Sigma \subset X$,
subject to the only condition that their complementary $X \setminus
\Sigma$, that is, the set of regular, or in other words,
non-singular points, be dense in $X$. For instance, if $X =
\mathbb{R}^n$ is an Euclidean space, then the set $\Sigma \subset X$
of singularities can be the set of all points with at least one
irrational coordinate. Indeed, in this case the set $X \setminus
\Sigma$ of non-singular, or regular points is the set of points with
all coordinates rational numbers, thus it is dense in $X$. A
relevant and rather remarkable fact to note in this case is that the
cardinal of the singularity set $\Sigma$ is  {\it strictly larger}
than the cardinal of the set of non-singular points, namely, $X
\setminus \Sigma$.

\item Second, there is no restriction on the behaviour of functions $f
: X \longrightarrow E$ in the neighbourhood of points in their
singularity sets $\Sigma \subset X$.

\end{itemize}

Related to this second freedom in dealing with singularities, one
should not forget its significant importance in applications. Indeed, as
stated in the mentioned Picard Great Theorem, an analytic function in the
neighbourhood of an isolated singularity point which is an essential
singularity takes on all possible complex values infinitely often,
with at most a single exception. \\ 
Consequently, in the neighbourhood
of a singularity, one can expect a rather arbitrary behaviour when
one deals with more general functions than
analytic ones. \\

The {\it first} class of algebras of generalized functions was aimed
to deal with singularities within a systematic and as widely
applicable as possible nonlinear theory, [12-29]. This class
contains as {\it particular cases} all the subsequent classes of
differential algebras of generalized functions constructed so far.
Within this largest class, a special subclass - of so called {\it nowhere
dense} algebras - was developed from the beginning, class which is able to
deal with arbitrary closed nowhere dense subsets $\Gamma \subset X$
of singularities, while no restrictions whatsoever are imposed on
the behaviour of generalized functions in the
neigbourhood of singularities. \\

The {\it second} class of algebras, [31-37], requires polynomial
growth conditions on generalized functions in the neigbourhood of
singularities. \\
In this regard, these algebras of generalized functions - which are
but a particular case of the infinite variety of all possible
differential algebras of generalized functions introduced in [12-29]
- suffer from a severe limitation. Namely, in the neighbourhood of
singularities of their generalized functions, these algebras require
a polynomial type growth condition, thus they cannot deal even with
isolated singularities such as essential singularities of analytic
functions. \\

The {\it third} class of algebras, [38-41, 9-11], is much more
powerful than the class of so called {\it nowhere dense} ones, since
these algebras are able to deal with arbitrary subsets $\Sigma
\subset X$ of singularities, subject to the mild condition that the
respective complementary subsets $X \setminus \Sigma$ be dense in
$X$, while again, no restrictions whatsoever are imposed on the
behaviour of generalized functions in the
neigbourhood of singularities in $\Sigma$. \\
An important fact to note here is that the subsets $\Sigma$ of
singularities can have a cardinal {\it larger} than that of the
subsets $X \setminus \Sigma$ of nonsingular or regular points, since
the condition that $X \setminus \Sigma$ be dense in $X$ can be
satisfied even when $X \setminus \Sigma$ is merely a dense countable
subset of $X$, in which case $\Sigma$
must of course be uncountable. \\

As for the nonlinear operations on singularities, the nowhere dense algebras and those in the third class allow arbitrary
smooth such operations, while in the second class only smooth operations with polynomial growth are possible. \\

As a consequence of its restriction upon singularities, as well as upon operations on singularities, the second class of
algebras cannot deal with a number of important problems which are easily treated within the nowhere dense algebras, or
those in the
third class. Among such problems are the following. \\

The global version of the classical Cauchy-Kovalevskaia theorem for solutions of analytic systems of nonlinear partial
differential equations cannot even be formulated, let alone solved, within the second class of algebras. \\
On the other hand, the first class of algebras is already able to produce such a global version on the existence of
solutions, [19,21,22]. \\

Arbitrary Lie group actions, which are of major importance in the solution of partial differential equations cannot be
defined within the second class of algebras. \\
Here again, the nowhere dense algebras are already enough to define globally arbitrary Lie group actions. And as one of
the consequences, one can for the first time obtain the complete solution of Hilbert's Fifth Problem, [28]. \\
This problem, again, cannot be formulated, let alone solved, within
the second class of algebras due to the polynomial mentioned
type growth conditions which they require. \\

Also, when defining differential algebras of generalized functions in the case of domains $X$ which are arbitrary finite dimensional smooth manifolds, the algebras in the first and third classes allow for considerably simpler constructions than those in the second class. \\

At a deeper analysis, however, one that is done in terms of Sheaf Theory, the essential difference between the nowhere
dense algebras or those in the third class, and on the other hand, the algebras in the second class, is that the former are
{\it flabby sheaves}, while the latter fail to be so, [9-11]. And as is known, [42], the lack of the flabbiness property
in the case of spaces of functions or generalized functions is an essential indicator of their limitations in dealing with
singularities. \\

Lastly, it should be noted that, in [43], the study of a fourth class which is far larger then the above third class of
algebras has been initiated. \\ \\

{\bf 1. Inclusion Diagrams and Reduced Power Algebras \\
        \hspace*{0.5cm} with the corresponding Ideals} \\

It is an elementary property of the linear vector space ${\cal D}' ( X )$ of Schwartz distributions that it can be represented as the {\it quotient} vector space \\

(1.1)~~~ $ {\cal D}' ( X ) = {\cal S}^\infty ( X ) / {\cal V}^\infty ( X ) $ \\

of the {\it vector subspaces} \\

(1.2)~~~ $ {\cal V}^\infty ( X ) ~\longrightarrow~ {\cal S}^\infty ( X ) ~\longrightarrow~
                                                                    ( {\cal C}^\infty ( X ) )^\mathbb{N} $ \\

with the arrows ''$\longrightarrow$'' representing usual inclusions ''$\subseteq$'', and \\

(1.3)~~~ $ {\cal S}^\infty ( X ) ~=~ \{~ s = ( \psi_\nu )_{\nu \in \mathbb{N}}
               \in ( {\cal C}^\infty ( X ) )^\mathbb{N} ~|~ s ~~\mbox{converges weakly in}~ {\cal D}' ( X ) ~\} $ \\

(1.4)~~~ $ {\cal V}^\infty ( X ) ~=~ \{~ v = ( \chi_\nu )_{\nu \in \mathbb{N}}
               \in {\cal S}^\infty ( X ) ~|~ v ~~\mbox{converges weakly to}~~ 0 ~~\mbox{in}~ {\cal D}' ( X ) ~\} $ \\

The remarkable fact, which has always been there in (1.1) - (1.4), is that in the right hand term of (1.2) we have the {\it differential algebra} $( {\cal C}^\infty ( X ) )^\mathbb{N}$, yet in (1.2) one uses only vector subspaces of it. Indeed, it takes little imagination to try to replace (1.2) with \\

(1.5)~~~ $ {\cal I}( X ) ~\longrightarrow~ {\cal A} ( X ) ~\longrightarrow~
                                                                    ( {\cal C}^\infty ( X ) )^\mathbb{N} $ \\

where ${\cal A} ( X )$ is a subalgebra in $( {\cal C}^\infty ( X ) )^\mathbb{N}$, while ${\cal I}( X )$ is an ideal in ${\cal A} ( X )$, and thus instead of the quotient vector space in (1.1), obtain the quotient {\it algebra} \\

(1.6)~~~ $ A ( X ) = {\cal A} ( X ) / {\cal I} ( X ) $ \\

which may allow a {\it nonlinear} theory of generalized functions, thus of singularities as well. \\

Indeed, for that purpose, it may be convenient to have the {\it inclusion}, that is, {\it linear embedding} \\

(1.7)~~~ $ {\cal D}' ( X ) ~\longrightarrow~ A ( X ) $ \\

and of course, also suitable partial derivations on $A ( X )$, which in some convenient manner may extend the distributional partial derivations on ${\cal D}' ( X )$. Clearly, such partial derivations can easily be obtained on $A ( X )$, in case ${\cal A} ( X )$ and ${\cal I} ( X )$ are {\it invariant} under the natural term-wise partial derivations \\

(1.8)~~~ $ ( {\cal C}^\infty ( X ) )^\mathbb{N} \ni s  = ( \psi_\nu )_{\nu \in \mathbb{N}} ~\longmapsto~
                         D^p s = ( D^p \psi_\nu )_{\nu \in \mathbb{N}} \in ( {\cal C}^\infty ( X ) )^\mathbb{N} $ \\

with $p \in \mathbb{N}^n$. Namely, if one has \\

(1.9)~~~ $ D^p {\cal I} ( X ) \subseteq {\cal I} ( X ),~~~~
                 D^p {\cal A} ( X ) \subseteq {\cal A} ( X ),~~~~ p \in \mathbb{N}^n $ \\

then one can simply define for $p \in \mathbb{N}^n$, the corresponding partial derivation \\

(1.10)~~~ $ A ( X ) \ni s + {\cal I} ( X ) ~\longmapsto~ D^p s + {\cal I} ( X ) \in  A ( X ) $ \\

Let us for the moment, however, deal only with the algebraic aspects of (1.5) - (1.7). An obvious immediate and simple way to obtain
(1.5) - (1.7) would be to construct {\it inclusion diagrams} of the form, [12-29] \\

\begin{math}
\setlength{\unitlength}{0.2cm}
\thicklines
\begin{picture}(60,18)

\put(8,15){${\cal I} ( X )$}
\put(14,15.5){\vector(1,0){16}}
\put(31,15){${\cal A} ( X )$}
\put(37,15.5){\vector(1,0){16}}
\put(54,15){$( {\cal C}^\infty ( X ) )^\mathbb{N}$}

\put(10,5){\vector(0,1){8}}
\put(32.5,5){\vector(0,1){8}}

\put(8,2){${\cal V}^\infty ( X )$}
\put(16,2.6){\vector(1,0){14}}
\put(31,2){${\cal S}^\infty ( X )$}

\put(0,8){$(1.11)$}

\end{picture}
\end{math} \\

which satisfy the condition \\

(1.12)~~~ $ {\cal I} ( X ) \bigcap {\cal S}^\infty ( X ) ~=~ {\cal V}^\infty ( X ) $ \\

a condition which is both {\it necessary} and {\it sufficient} for the existence of the {\it linear embedding} \\

(1.13)~~~ $ {\cal S}^\infty ( X ) / {\cal V}^\infty ( X ) \ni s + {\cal V}^\infty ( X ) ~\longmapsto~
                  s + {\cal I} ( X ) \in {\cal A} ( X ) / {\cal I} ( X ) $ \\

thus equivalently, for (1.7). \\

Unfortunately however, inclusion diagrams (1.11) {\it cannot} be constructed in view of the simple fact that, [12-29] \\

(1.14)~~~ $ (~ {\cal V}^\infty ( X ) ~.~ {\cal V}^\infty ( X ) ~) \bigcap {\cal S}^\infty ( X )
                                                                \nsubseteqq {\cal V}^\infty ( X ) $ \\

as simple counterexamples can show it. Indeed, it is easy to construct sequences $v = ( \chi_\nu )_{\nu \in \mathbb{N}} \in {\cal V}^\infty ( X )$, such that $v^2 \in {\cal S}^\infty ( X )$, yet $v^2 \notin {\cal V}^\infty ( X )$. For instance, when $X = \mathbb{R}$, one can take $\chi_\nu ( x ) = \cos ( \nu x)$, and obtain indeed that $v \in {\cal V}^\infty ( X )$, $v^2 \in {\cal S}^\infty ( X )$, and furthermore $v^2$ converges weakly to $1 / 2$ in ${\cal D}' ( X )$, thus clearly $v^2 \notin {\cal V}^\infty ( X )$. \\

Consequently, one can turn to the immediately more involved inclusion diagrams, [12-29] \\

\begin{math}
\setlength{\unitlength}{0.2cm}
\thicklines
\begin{picture}(60,31)

\put(8,28){${\cal I} ( X )$}
\put(14,28.5){\vector(1,0){16}}
\put(31,28){${\cal A} ( X )$}
\put(37,28.5){\vector(1,0){16}}
\put(54,28){$( {\cal C}^\infty ( X ) )^\mathbb{N}$}

\put(10,18){\vector(0,1){8}}
\put(32.5,18){\vector(0,1){8}}

\put(9,15){${\cal V}$}
\put(12,15.5){\vector(1,0){18}}
\put(32,15){${\cal S}$}

\put(10,13.5){\vector(0,-1){8}}
\put(32.5,13.5){\vector(0,-1){8}}

\put(8,2){${\cal V}^\infty ( X )$}
\put(16,2.6){\vector(1,0){14}}
\put(31,2){${\cal S}^\infty ( X )$}

\put(0,15){$(1.15)$}

\end{picture}
\end{math} \\

where ${\cal V}$ and ${\cal S}$ are vector subspaces, such that the following three conditions hold \\

(1.16)~~~ $ {\cal I} ( X ) \bigcap {\cal S} ~=~ {\cal V} $ \\

(1.17)~~~ $ {\cal V}^\infty ( X ) \bigcap {\cal S} ~=~ {\cal V} $ \\

(1.18)~~~ $ {\cal V}^\infty ( X ) + {\cal S} ~=~ {\cal S}^\infty ( X ) $ \\

which, as it is easy to see, are both {\it necessary} and {\it sufficient} for the existence of the {\it linear embeddings} \\

(1.19)~~~ $ {\cal D}' ( X ) S ~\longmapsto~ s + {\cal V}^\infty ( X ) \in
                            {\cal S}^\infty ( X ) / {\cal V}^\infty ( X ) $ \\

(1.20)~~~ $ {\cal S} / {\cal V} \ni s + {\cal V} ~\longmapsto~ s + {\cal V}^\infty ( X ) \in
                            {\cal S}^\infty ( X ) / {\cal V}^\infty ( X ) $ \\

(1.21)~~~ $ {\cal S} / {\cal V} \ni s + {\cal V} ~\longmapsto~
                              s + {\cal I} ( X ) \in A ( X ) = {\cal A} ( X ) / {\cal I} ( X ) $ \\

where the mappings (1.19), (1.20) are in fact {\it vector space isomorphisms}. \\

Now clearly, (1.19) - (1.21) give the desired {\it linear embedding} (1.7) of the Schwartz distributions into algebras of generalized functions, namely \\

(1.22)~~~ $ {\cal D}' ( X ) ~\longrightarrow~ A ( X ) ~=~ {\cal A} ( X ) / {\cal I} ( X ) $ \\

In view of (1.5), (1.6), the algebras of generalized functions $A ( X )$ in (1.7), (1.22) are nothing else but {\it
reduced powers} of the algebra ${\cal C}^\infty ( X )$ of smooth functions on $X$. \\

The general Model Theoretic, [44], construction of reduced powers, although hardly known as such among so called working
mathematicians, happens nevertheless to appear in quite a number of important places in Mathematics at large. For a sample
of them, one can note the following. The Cauchy-Bolzano construction of the field $\mathbb{R}$ of usual real numbers is in
fact a reduced power of the rational numbers $\mathbb{Q}$. More generally, the completion of any metric space is a reduced
power of that space. Furthermore, this is but a particular case of the fact that the completion of any uniform topological
space is a reduced power of that space. Also, in a rather different direction, the field $^*\mathbb{R}$ of nonstandard
real numbers can be obtained as a reduced power of the usual field $\mathbb{R}$ of real numbers. \\

In view of the above, the use of reduced powers in the construction of algebras of generalized functions should not be seen as much else but a further application of that basic construction in Model Theory, this time to the study of large classes of singularities, and as such, to the solution of rather general nonlinear systems of partial differential equations, among others. \\

As seen in the sequel, the {\it ideals} ${\cal I} ( X )$ in (1.22) play the essential role in the inevitable infinite branching which occurs when multiplying generalized functions that have singularities above a certain level. \\ \\

{\bf 2. Off-Diagonality Characterization} \\

A fundamental result in the nonlinear algebraic theory of
generalized functions, see 46F30, is the simple and purely algebraic
{\it characterization} of the algebras of generalized functions
(1.7) which are constructed upon inclusion diagrams (1.15) - (1.21).
In this regard, first we note that these inclusion diagrams can be
simplified as follows. In the inclusion diagrams (1.15) there are
{\it four} spaces to be chosen, namely ${\cal I} ( X ), {\cal A} ( X
), {\cal V}$ and ${\cal S}$. However, it is easy to see, [12-29],
that such inclusion diagrams can be reduced to the simpler form
depending
only on {\it two} spaces $( {\cal V}, {\cal S}\,' )$, namely \\

\begin{math}
\setlength{\unitlength}{0.2cm}
\thicklines
\begin{picture}(60,31)

\put(6,28){${\cal I} ( {\cal V}, {\cal S}\,' )$}
\put(15,28.5){\vector(1,0){13}}
\put(29,28){${\cal A} ( {\cal V}, {\cal S}\,' )$}
\put(38,28.5){\vector(1,0){15}}
\put(54,28){$( {\cal C}^\infty ( X ) )^\mathbb{N}$}

\put(10,18){\vector(0,1){8}}
\put(32.5,18){\vector(0,1){8}}

\put(9,15){${\cal V}$}
\put(12,15.5){\vector(1,0){17}}
\put(31,15){${\cal V} \oplus {\cal S}\,'$}
\put(53,15.5){\vector(-1,0){15}}
\put(55,15){${\cal U}^\infty ( X )$}

\put(10,13.5){\vector(0,-1){8}}
\put(32.5,13.5){\vector(0,-1){8}}

\put(8,2){${\cal V}^\infty ( X )$}
\put(16,2.6){\vector(1,0){14}}
\put(31,2){${\cal S}^\infty ( X )$}

\put(0,15){$(2.1)$}

\end{picture}
\end{math} \\

where ${\cal V}, {\cal S}\,'$ are vector subspaces in ${\cal S}^\infty ( X )$, while \\

(2.2)~~~ $ {\cal U}^\infty ( X ) ~=~
             \{~ u_\psi = ( \psi, \psi, \psi, \ldots ) ~|~ \psi \in {\cal C}^\infty ( X ) ~\}
                                          \subset {\cal S}^\infty ( X ) $ \\

is the {\it diagonal} in the Cartesian product $( {\cal C}^\infty ( X ) )^\mathbb{N}$. \\

As for the conditions (1-16) - (1.18), they now become \\

(2.3)~~~ $ {\cal V} \bigcap {\cal S}\,' ~=~ \{ 0 \} $ \\

(2.4)~~~ $ {\cal I} ( {\cal V}, {\cal S}\,' ) \bigcap {\cal S}\,' ~=~ \{ 0 \} $ \\

(2.5)~~~ $ {\cal V}^\infty ( X ) \oplus {\cal S}\,' ~=~ {\cal S}^\infty ( X ) $ \\

The mentioned fundamental result regarding the {\it characterization} of algebras of generalized functions (1.7)
constructed upon inclusion diagrams (1.15) - (1.21) as simplified in (2.1) - (2.5) is the following : \\

{\bf Theorem 2.1.} [12-29] \\

Within a large class of inclusion diagrams (2.1) - (2.5), the conditions (2.3) - (2.5) are equivalent with the following
{\it off-diagonality} one \\

(2.6)~~~ $ {\cal I} ( {\cal V}, {\cal S}\,' ) \bigcap {\cal U}^\infty ( X ) ~=~ \{ 0 \} $ \\

{\bf Remark 2.1.} \\

It is both theoretically and practically important to note that, as
seen in [12-29], there are {\it infinitely} many different inclusion
diagrams (2.1) - (2-5). Moreover, they give {\it infinitely} many
different corresponding algebras
of generalized functions \\

(2.7)~~~ $ A ~=~ {\cal A} ( {\cal V}, {\cal S}\,' ) / {\cal I} ( {\cal V}, {\cal S}\,' ) $ \\

which, see (1.22), contain the vector space ${\cal D}' ( X )$ of Schwartz distributions. \\ \\

{\bf 3. Inevitable Infinite Branching in the \\
        \hspace*{0.45cm} Multiplication of Singularities} \\

And now, we can come to the main issue in this paper, namely, to indicate the reason for the {\it inevitable} infinite
possibilities in defining multiplication of generalized functions in case their singularities are above certain levels. \\

The remarkable fact in this regard is that the respective reason is of a simple {\it algebraic} nature, namely, it is a
direct consequence of the off-diagonality characterization in Theorem 2.1. above of the algebras of generalized functions
(2.7 )constructed through the method of reduced powers. \\

Indeed, for a given pair ${\cal V}, {\cal S}\,'$ and a corresponding subalgebra ${\cal A} ( {\cal V}, {\cal S}\,' )
\subseteq ( {\cal C}^\infty ( X ) )^\mathbb{N}$ in an inclusion diagram (2.1) - (2-5), let us denote by \\

(3.1)~~~ $ \textbf{ID} ( X, {\cal V}, {\cal S}\,', {\cal A} ( {\cal V}, {\cal S}\,' ) ) $ \\

the set of all ideals ${\cal I}$ in ${\cal A} ( {\cal V}, {\cal
S}\,' )$ which can occur in such inclusion diagrams (2.1) - (2-5).
This means, therefore, that for every such ideal ${\cal I} \in
\textbf{ID} ( X, {\cal V}, {\cal S}\,',
{\cal A} ( {\cal V}, {\cal S}\,' ) )$, there exists a corresponding algebra of generalized functions \\

(3.2)~~~ $ A = {\cal A} ( {\cal V}, {\cal S}\,' ) / {\cal I} $ \\

which, in view of (1.22), contains the vector space ${\cal D}' ( X )$ of Schwartz distributions. \\

Now, the essential point regarding the multiplication of singularities is that, in view of the inclusion \\

(3.3)~~~ $ {\cal V} \oplus {\cal S}\,' \supset {\cal U}^\infty ( X ) $ \\

in (2.1), it follows easily that the multiplication in each of the algebras of generalized functions $A$ in (3.2)
preserves the usual multiplication of smooth functions in ${\cal C}^\infty ( X )$. \\

On the other hand, regarding the multiplication of generalized functions that are not smooth - therefore, are elements in
$A \setminus {\cal C}^\infty ( X )$ - it is well known that in general they do no longer preserve even the usual
multiplication of continuous functions, this being one of the immediate implications of the 1954 so called Schwartz
impossibility result, [12-29]. \\

Furthermore, in view of Remark 2.1. above, there are {\it
infinitely} many ways according to which multiplication ends up
being done, ways corresponding to the various algebras $A$ of
generalized functions. And the possibility of this infinite
branching of multiplication is manifested as soon as the generalized
functions which
are the factors in multiplication belong to $A \setminus {\cal C}^\infty ( X )$, and as such, and are farther and farther
removed from ${\cal C}^\infty ( X )$, or even form ${\cal C}^0 ( X )$, that is, have a {\it higher levels} of
singularity. \\

Let us illustrate the above with a simple example. For that purpose,
let us fix the pair ${\cal V}, {\cal S}\,'$ and take as a
corresponding subalgebra ${\cal A} ( {\cal V}, {\cal S}\,' ) = (
{\cal C}^\infty ( X ) )^\mathbb{N}$ in the inclusion
diagram (2.1) - (2-5). \\
In such a case, one may expect a natural or canonical multiplication
if one would be able to {\it single out} in a suitable manner a
certain ideal ${\cal I} \in \textbf{ID} ( X, {\cal V}, {\cal S}\,',
( {\cal C}^\infty ( X ) )^\mathbb{N} )$, and thus obtain the corresponding algebra of generalized functions $A$ in
(3.2). \\

Here however, the off-diagonality condition (2.6) interferes, leading to a rather involved structure for the set of ideals
which satisfy that condition, as seen in \\

{\bf Proposition 3.1.} \\

There is {\it no} largest ideal in the set \\

(3.4)~~~ $ \textbf{ID} ( X ) $ \\

of ideals ${\cal I}$ in $( {\cal C}^\infty ( X ) )^\mathbb{N}$ which satisfy the off-diagonality condition \\

(3.5)~~~ $ {\cal I} \bigcap {\cal U}^\infty ( X ) ~=~ \{ 0 \} $ \\

{\bf Proof.} \\

Let us again take $X = \mathbb{R}$, together with $v\,' = ( \chi\,'_\nu )_{\nu \in \mathbb{N}}, v\,''  =
( \chi\,''_\nu )_{\nu \in \mathbb{N}} \in ( {\cal C}^\infty ( X ) )^\mathbb{N}$, where \\

$~~~~~~ \chi\,'_\nu ( x ) = 1 + \sin ( \nu x ),~~~ \chi\,''_\nu ( x ) = 1 + \cos ( \nu x ),~~~
                                                                    \nu \in \mathbb{N}, x \in X $ \\

Then it follows easily that \\

$~~~~~~ {\cal I}\,' = v\,' ( {\cal C}^\infty ( X ) )^\mathbb{N},~
                      {\cal I}\,'' = v\,'' ( {\cal C}^\infty ( X ) )^\mathbb{N} \in \textbf{ID} ( X ) $ \\

However \\

$~~~~~~ {\cal I} ~=~ {\cal I}\,' + {\cal I}\,'' $ \\

is an ideal in $( {\cal C}^\infty ( X ) )^\mathbb{N}$ which fails to satisfy the off-diagonality condition (3.5),
therefore \\

(3.6)~~~ $ {\cal I} \notin \textbf{ID} ( X )  $ \\

Indeed

$~~~~~~ v = v\,' + v\,'' \in {\cal I} $ \\

and $v = ( \psi_\nu )_{\nu \in \mathbb{N}}$, where \\

$~~~~~~ \psi_\nu ( x ) = 2 + \sin ( \nu x ) + \cos ( \nu x ) > 0,~~~ \nu \in \mathbb{N}, x \in X $ \\

Consequently \\

$~~~~~~ {\cal I} ~=~ ( {\cal C}^\infty ( X ) )^\mathbb{N} $ \\

thus it is not a proper ideal in $( {\cal C}^\infty ( X ) )^\mathbb{N}$, and in particular, it does not satisfy condition
(3.5). \\

It follows that $\textbf{ID} ( X )$ does not contain ideals which may contain both ideals ${\cal I}\,'$ and
${\cal I}\,''$.

\hfill $\Box$ \\

In view of [12-29], it is obvious that there are infinitely many
ideals in  $\textbf{ID} ( X )$, and in fact, also in ${\cal I} \in
\textbf{ID} ( X, {\cal V}, {\cal S}\,', ( {\cal C}^\infty ( X )
)^\mathbb{N} )$. Furthermore, various ideals in these sets clearly
lead to significantly different multiplications in the corresponding
algebras of generalized functions (3.2). \\ \\

{\bf 4. Adhock Multiplications} \\

There is a considerable literature, [34], in which multiplications,
and some more general nonlinear operations are performed upon a few
specific generalized functions with singularities, often including
the Dirac delta distribution and possibly, its derivatives. Such
studies, mostly undertaken by, and of interest to certain physicists
and engineers, do not assume, let alone, aim to, or develop any
general nonlinear theory for dealing with singularities, and
instead, limit themselves clearly to
a few special instances of such operations. \\
The respective authors seem not to be cognisant of the wider existing
relevant literature, let alone of the decades long existing systematic
nonlinear theory of generalized functions. \\ \\

{\bf 5. The Apparent Novelty of the Inevitable Infinite \\
        \hspace*{0.5cm} Branching in the Multiplication of Singularities} \\

The remarkable survey of state of the art Physics, and also of a lot
of the Mathematics underlying it in [45], see review [46], mentions
in passing the Schwartz distributions. However, in its declared
intention of guide to the laws of the universe, fails to point to
the issue of nonlinear operations on singularities, let alone, to
mention the inevitable infinite branching which appears when
multiplication or more general nonlinear operations on singularities
are
performed. \\

On the other hand, most important physical phenomena and processes
lead to singularities in their mathematical formulation.
Furthermore, the respective mathematical models happen essentially
to be nonlinear. One of them is the phenomenon of black holes in
General Relativity, a phenomenon to which the author of [45] has
contributed in significant ways. \\

And yet, when it comes to dealing with singularities in highly
nonlinear contexts, there seems not to be a sufficient awareness
about such critically important novel features as the inevitable
infinite branchings which appear in such situations, when
multiplication or more general nonlinear operations
on singularities are performed. \\


\begin{thebibliography}{99}

\bibitem{} Schwartz L : Sur l'impossibilite de la multiplication des distributions.
C.R. Acad. Sci. Paris, vol. 239, 1954, 847-848

\bibitem{} Hoermander L : Linear Partial Differential Operators. (fourth printing) Springer, New York, 1976

\bibitem{} Mallios A : Geometry of Vector Sheaves. An Axiomatic Approach to
Differential Geometry, vols. I (chapts. 1-5), II (chapts. 6-11)
Kluwer, Amsterdam, 1998

\bibitem{} Mallios A : Modern Differential Geometry in Gauge Theories.
Volume 1 : Maxwell Fields, Volume 2 : Yang-Mills Fields. Birkhauser,
Boston, 2006

\bibitem{} Mallios A : On an axiomatic treatment of differential geometry
via vector.sheaves. Applications. (International Plaza) Math. Japon-
ica, vol. 48, 1998, pp. 93-184

\bibitem{} Mallios A : The de Rham-Kahler complex of the Gelfand sheaf of a
topological algebra. J. Math. Anal. Appl., 175, 1993, 143-168

\bibitem{} Mallios A : On an abstract form of Weils intergrality theorem.
Note Mat., 12, 1992, 167-202 (invited paper)

\bibitem{} Mallios A : On an axiomatic approach to geometric prequanti-
zation : A classification scheme a la Kostant-Souriau-Kirillov.
J.Mathematical Sciences (former J. Soviet Math.), vol. 98-99

\bibitem{} Mallios A, Rosinger E E : Abstract differential geometry, differential algebras of generalized functions, and de Rham cohomology. Acta Applicandae Mathematicae, vol. 55, no. 3, February 1999, 231 – 250

\bibitem{} Mallios A, Rosinger E E : Space-time foam dense singularities and de Rham cohomology. Acta Applicandae Mathematicae, vol. 67, no. 1, May 2001, 59 – 89

\bibitem{} Mallios A, Rosinger E E : Dense singularities and de Rham cohomology. In (Eds. Strantzalos P, Fragoulopoulou M) Topological Algebras with Applications to Differential Geometry and Mathematical Physics. Proc. Fest-Colloq. in honour of Prof. Anastasios Mallios (16-18 September 1999), pp. 54-71, Dept. Math. Univ. Athens Publishers, Athens, Greece, 2002

\bibitem{} Rosinger E E : Embedding of the D' distributions into pseudotopological algebras. Stud. Cerc. Math., Vol. 18, No. 5, 1966, 687-729

\bibitem{} Rosinger E E : Pseudotopological spaces, the embedding of the D' distributions into algebras. Stud. Cerc. Math., Vol. 20, No. 4, 1968, 553-582

\bibitem{} Rosinger E E : Pseudotopological Structures. A Uniform Continuous Embedding of the D' Distributions into Pseudotopological Algebras. Dr.Sc. Thesis, University of Bucharest, Romania, 1972. Advisor Professor G. Marinescu

\bibitem{} Rosinger E E : Division of Distributions. Pacif.J. Math., Vol. 66, No. 1, 1976, 257-263

\bibitem{} Rosinger E E : Nonsymmetric Dirac distributions in scattering theory. In Springer Lecture Notes in Mathematics, Vol. 564, 391-399, Springer, New York, 1976

\bibitem{} Rosinger E E : Distributions and Nonlinear Partial Differential Equations. Springer Lectures Notes in Mathematics, Vol. 684, Springer, New York, 1978

\bibitem{} Rosinger E E : Nonlinear Partial Differential Equations, Sequential and Weak Solutions, North Holland Mathematics Studies, Vol. 44, Amsterdam, 1980

\bibitem{} Rosinger E E : Generalized Solutions of Nonlinear Partial Differential Equations. North Holland Mathematics Studies, Vol. 146, Amsterdam, 1987

\bibitem{} Rosinger E E : An exact computation of the Feynman path integrals. Technical Report, 1989, Department of Mathematics and Applied Mathematics, University of Pretoria

\bibitem{} Rosinger E E : Nonlinear Partial Differential Equations, An Algebraic View of Generalized Solutions. North Holland Mathematics Studies, Vol. 164, Amsterdam, 1990

\bibitem{} Rosinger E E : Global version of the Cauchy-Kovalevskaia theorem for nonlinear PDEs. Acta Appl. Math., Vol. 21, 1990,

\bibitem{} Rosinger E E : Characterization for the solvability of nonlinear PDEs, Trans. AMS, Vol. 330, 1992

see also reviews MR 92d:46098, Zbl. Math. 717 35001, MR 92d:46097, Bull. AMS vol.20, no.1, Jan 1989, 96-101, MR
89g:35001

\bibitem{} Rosinger E E, Walus E Y : Group invariance of generalized solutions obtained through the algebraic method. Nonlinearity, Vol. 7, 1994, 837-859

\bibitem{} Rosinger E E, Walus E Y : Group invariance of global generalised solutions of nonlinear PDEs in nowhere dense algebras. Lie Groups and their Applications, Vol. 1, No. 1, July-August 1994, 216-225

\bibitem{} Rosinger E E : Nonprojectable Lie Symmetries of nonlinear PDEs and a negative solution to Hilberts fifth problem. In (Eds. N.H. Ibragimov and F.M. Mahomed) Modern Group Analysis VI, Proceedings of the International Conference in the New South Africa, Johannesburg, January 1996, 21-30. New Age Inter. Publ., New Delhi, 1997

\bibitem{} Rosinger E E : Parametric Lie Group Actions on Global Generalised Solutions of Nonlinear PDEs, Including a Solution to Hilberts Fifth Problem. Kluwer Acad. Publishers, Dordrecht, Boston, 1998

\bibitem{} Rosinger E E : Arbitrary Global Lie Group Actions on Generalized Solutions of Nonlinear PDEs and an Answer to Hilberts Fifth Problem. In (Eds. Grosser M, H¨ormann G, Kunzinger M, Oberguggenberger M B) Nonlinear Theory of Generalized Functions, 251-265, Research Notes in Mathematics, Chapman \& Hall/CRC, London, New York, 1999

\bibitem{} Rosinger E E : Scattering in highly singular potentials. \\
arXiv:quant-ph/0405172

\bibitem{} Rosinger E E : Junction Conditions, Resolution of Singularities and Nonlinear Equations of Physics.
arXiv:math/0611445

\bibitem{} Colombeau J-F : New Generalized Functions and Multiplication
of Distributions. Mathematics Studies, Vol. 84, North-Holland,
Amsterdam, 1984

\bibitem{} Colombeau J-F : Elementary Introduction to New Generalized
Functions. North-Holland Mathematics Studies, Vol. 113, Amsterdam,
1985

\bibitem{} Biagioni H A : A Nonlinear Theory of Generalized Functions.
Lecture Notes in Mathematics, Vol. 1421, Springer, New York, 1990

\bibitem{} Oberguggenberger M B : Multiplication of Distributions and Applications to PDEs. Pitman Research Notes in Mathematics, Vol. 259, Longman, Harlow, 1992

\bibitem{} Grosser M, Hoermann G, Kunzinger M, Oberguggenberger M B
(Eds.) : Nonlinear Theory of Generalized Functions. Pitman Research
Notes in Mathematics, Chapman \& Harlow / CRC , London, 1999

\bibitem{} Grosser M, Kunzinger M, Oberguggenberger M, Steinbauer R :
Geometric Theory of Generalized Functions with Applications to
General Relativity. Kluwer, Dordrecht, 2002

\bibitem{} Colombeau J F : Mathematical problems on generalized functions
and the canonical Hamiltonian formalism. arXiv:0708.3425

\bibitem{} Rosinger E E : How to solve smooth nonlinear PDEs in algebras of generalized functions with dense singularities (invited paper) Applicable Analysis, vol. 78, 2001, 355-378

\bibitem{} Rosinger E E : Differential Algebras with Dense Singularities on Manifolds. Acta Applicandae Mathematicae.
Vol. 95, No. 3, Feb. 2007, 233-256, arXiv:math/0606358

\bibitem{} Rosinger E E : Which are the Maximal Ideals ? \\ arXiv:math.GM/0607082

\bibitem{} Rosinger E E : Space-Time Foam Differential Algebras of Generalized Functions and a Global Cauchy-Kovalevskaia
Theorem. Acta Applicandae Mathematicae, DOI
10.1007/s10440-008-9326-z, Received: 5 February 2008 / Accepted : 23
September 2008

\bibitem{} Kaneko A : Introduction to Hyperfunctions. Kluwer, Dordrecht, 1088

\bibitem{} Rosinger E E : Local Functions : Algebras, Ideals, and Reduced
Power Algebras. arXiv:0912.4049

\bibitem{} Bell J L, Slomson A B : Models and Ultraproducts, An Introduction.
North-Holland, Amsterdam, 1969

\bibitem{} Penrose R : The Road to Reality, a Complete Guide to the Laws
of the Universe. Jonathan Cape, London, 2004.

\bibitem{} Rosinger E E : The Road to Reality. Letter to the Editor, Mathematical
Intelligencer, Vol. 29, No. 3, 2007, pp. 5,6

\bibitem{} Synowiec J A : Some highlights in the development of Algebraic Analysis.
Algebraic Analysis and Related Topics. Banach Center Publications,
Vol. 53, 2000, 11-46, Polish Academy of Sciences, Warszawa

\bibitem{} Lewy, H : An example of smooth linear partial differential equation
without solutions. Ann. Math., vol. 66, no. 2, 1957, 155-158

\end{thebibliography}
\end{document}